\documentclass[a4paper,12pt]{article}
\usepackage[T2A]{fontenc}
\usepackage[cp1251]{inputenc}
\usepackage[english,russian]{babel}
\usepackage{graphicx}
\usepackage{amsmath}
\usepackage{color}
\usepackage{amssymb}
\usepackage{chicago}
\usepackage{graphicx}
\usepackage{epstopdf}
\usepackage{epsfig}
\usepackage{multirow}

\newtheorem{theorema}{Теорема}

\newtheorem{ddef}{Определение}

\newtheorem{lemma}{Лемма}

\newenvironment{proof}{\par\noindent{\bf Доказательство.}}{\hfill$\scriptstyle\blacksquare$}
\makeatletter
\renewcommand \thesection {\@arabic\c@section.}
\renewcommand\thesubsection {\thesection\@arabic\c@subsection.}
\renewcommand\thesubsubsection {\thesubsection\@arabic\c@subsubsection.}
\renewcommand\theparagraph {\thesubsubsection\@arabic\c@paragraph.}

\renewcommand{\ge}{\geqslant}
\renewcommand{\le}{\leqslant}

\renewcommand{\kappa}{\varkappa}
\renewcommand{\epsilon}{\varepsilon}
\renewcommand{\phi}{\varphi}

\begin{document}
\selectlanguage{russian}

\title{
{\small \begin{flushright}
УДК 519.714
\end{flushright}}
Сложность дизъюнктивных нормальных форм и полу-эффект Шеннона в некоторых подклассах булевых функций.
}

\author{С.~С. Гранин \\ МФТИ \and Ю.~В. Максимов\thanks{Исследование выполнено при финансовой поддержке РФФИ в рамках научного проекта \No\,14-07-31277 мол\_а; а также при частичной поддержке Лаборатории структурных методов анализа данных в предсказательном моделировании ФУПМ МФТИ, грант правительства РФ дог. 11.G34.31.0073.} \\ ПреМоЛаб МФТИ и ИППИ РАН
}
\date{}

\maketitle

\begin{abstract}
В работе рассматривается задача построения минимальных и кратчайших дизъюнктивных нормальных форм булевых функций с ограниченным числом. В случае полиномиально ограниченного множества нулей получены новые рекордные нижние оценки сложности почти всех функций соответствующих классов. Указанные оценки и в сочетании с известными раннее верхними оценками позволяют установить точный порядок роста сложности ДНФ булевых функций с не более чем полиномиальным числом нулей. В работе также ставиться вопрос об эффекте Шеннона в соответствующих малому числу подклассах функций.

\vspace{0.5cm} \noindent \textbf{Ключевые слова:}
задача покрытия множеств, булева функция, дизъюнктивная нормальная форма, схемная сложность.
\end{abstract}

\section{Введение}
При решении ряда современных задач дискретной математики и информатики, таких как построение эффективных $SAT$-солверов или решения $NP$-трудных задач о покрытии множеств, часто возникает задача решения булевых уравнений Нельсоновского типа, представимых в виде
\begin{gather}\label{eq:nelson}
\prod\limits_{i=1}^q F_i = 0,
\end{gather}
где каждая из формул $F_i$ над пространством булевых переменных $x_1,x_2,\dots$, $x_n$ записана в виде $F_i = \bigvee\limits_{i\in I} x_{i}^{\sigma_i}$ для некоторого $I\subseteq [n]$.

Одним из наиболее интересных случаев указанной постановки является случай, в котором существенная часть функций $F_i$ отвечает функциям обращающимся в ноль в единственной точке.

В общем случае указанная задача сводиться к так называемой задаче исключения нулей из дизъюнктивной нормальной формы (ДНФ), отвечающей частичному произведению функций $F_i$, $i\in I$.
В настоящей работе будут получены нижние оценки и сложности функций, превосходящие, в определенных случаях, известные ранее, в частности  \cite{dj01,dj02}, а также \cite{mu06,zk85dan,zk86}, \cite{kog87,me12a,me12b,me13} и другие.

\section{Структура работы}
В первой части работы рассматривается общий метод построения высоких нижних оценок дизъюнктивных нормальных форм булевых функций, заложенный в работе \cite{dj01} и развитый в работах \cite{me12a,me12b,me13}.

Суть метода состоит в геометрическом анализе исходной булевой функции. А именно, проблема построения простой ДНФ сводиться к задаче построения покрытия не всего множества точек, отвечающих единицам рассматриваемой функции, а некоторого его подмножества малой мощности. В случае, если удается выделить подмножество для покрытия которого, при условии его небольшой мощности, необходимо достаточно большое число конъюнкций, то указанный метод позволяет получать нижние оценки сложности существенно превосходящие имеющиеся методы.

В следующей части нами приводиться обзор известных результатов, касающихся поведения типичных и экстремальных по сложности функций, в зависимости от числа точек, на которых они обращаются в ноль (см. таблицу \ref{table:review}).

В заключительной части работы нами представлен список теоретических вопросов, касающихся задачи построения минимальных ДНФ функций и смежных вопросов, ответы на которые представляются авторам интересными и важными.

\section{Терминология и обозначения}

Предполагается, что читатель знаком с основными понятиями теории дизъюнктивных нормальных форм и теории вероятности. Основную часть неопределенных здесь понятий, относящиеся к теории построения дизъюнктивных нормальных форм, определены в работе \cite{ja10}.

Задача минимизации сложности дизъюнктивной нормальной формы булевой функции, заданной в конъюнктивной нормальной форме (КНФ) уравнением \ref{eq:nelson}, состоит в том, что бы записать указанную функцию в дизъюнктивной нормальной форме (ДНФ) с использованием как можно меньшего числа логических произведений (конъюнкций), а также символов переменных и их отрицаний (литералов). \emph{Длиной} ДНФ называется число входящих в нее конъюнкций, \emph{рангом} ДНФ называется число использованных в её записи литералов.

Дизъюнктивная нормальная форма минимального ранга называется \emph{минимальной}, а минимальной длины --- \emph{кратчайшей}. С точки зрения построения оптимальных покрытий, ДНФ минимального длины соответствует покрытие использующее минимальное число граней; ДНФ минимального ранга соответствует минимальное вес совокупности граней, где в качестве веса грани выступает разность между размерностью функции и размерностью грани.

Литералом будем называть булеву переменную или ее отрицание. Литерал $x_i$ булевой функции $f(x_1,x_2,...,x_n)$ ассоциируем со столбцом $M_{[k]}^i$ матрицы нулей функции $f$; литерал $\bar{x_i}$ ассоциируем с покоординатным отрицанием столбца $M_{[k]}^i, 1\le i \le n$. Далее под термином литерал $x_i^{\sigma_i}$ будет часто подразумеваться ассоциированный с ним булев вектор.

Обозначим через $P_k^n$ класс булевых функций от $n$ переменных, имеющих ровно $k$ нулевых точек.

Обозначим через $\log x$ натуральный алгоритм числа $x$, $\log^{(1)} x = \log x$, $\log^{(k+1)} x = \log^{(k)} x$ при всех $k \ge 1$; асимптотические нотации $O(.), o(.)$ и $\omega(.)$ всегда, кроме специально оговоренных случаев, применяются при $n\rightarrow\infty$. Выражения $f \ll g$ и $f \gg g$ применяемые в неравенствах и формулировках утверждений, равносильны равенствам $f = o(g)$ и $f = \omega(g)$ соответственно.


\textit{Матрицей нулей} $M_f$ функции $f$ назовем булеву матрицу, по строкам которой расположены нули функции $f$. Обозначим $[k] = \{1,2,...,k\}$; $M_{i_1,i_2,...,i_t}^{j_1,j_2,...,j_r}$ подматрицу матрицы $M$ образованную пересечением строк $i_1$, $i_2$, ..., $i_t$ со столбцами $j_1, j_2, ...j_r$.

Скажем, что для \emph{почти всех булевых функций класса $P^n_{k(n)}$ выполнено свойство $A$}, если доля функций для которых оно не справдливо стремиться к 0, при $n\rightarrow\infty$.

Положим $B_k$ --- множество всех булевых векторов размерности $k$; $B_k^0$ и $B_k^1$ --- множество булевых векторов, первая координата которых равна 0 или 1 соответственно; под вектором $\bar{x}$ будем понимать покоординатное отрицание вектора $x\in B_k$; минимум из числа единиц и числа нулей $x$ назовем \textit{весом} вектора $x, x\in B_k$.

\section{Известные результаты}
Первые эффективные алгоритмы построения ДНФ булевых функций с ограниченным числом нулей были предложены в работах \cite{zk85dan,zk86}.

Предложенный им метод позволяет строить для почти всех функций из класса $P_k^n$ дизъюнктивные нормальные формы число конъюнкций которых, не превосходит $O\left(\frac{nk}{\log_2 n}(1 + o(1))\right)$. Указанные оценки во многих случаях остаются рекордными в смысле порядка роста функции.

А.~Ю. Коганом в работе \cite{kog87} были предложены способы получения нижних оценок сложности основанные на специальной релаксации исходной задачи к задаче линейного программирования. Ключевая для последующих рассуждений идея была представлена А.Г.~Дьяконовым в работе \cite{dj01}. Суть ее состоит в явном построении системы точек, покрытие которых трудно в классе ДНФ.

Основываясь в целом на методе А.Г. Дьяконова в ряде последующих работ были получены оценки для сложности булевых функций в классе ДНФ, являющихся в настоящий момент рекордными \cite{me12a,me12b,me13}.

\begin{table}[t]
\centering
\begin{tabular}{|c|l|l|}
  \hline
  Оценка & Автор & Комментарий \\
  \hline \hline
  $O(nk)$ & \cite{dj02}, & верхняя, произвольная функция\\
  & \cite{mu06}&\\
  \hline
  $O\left(\frac{nk}{\log_2 n}\right)$ & \cite{zk85dan} & верхняя, почти все функции, $k\le 2^{n/2}$ \\
  &\cite{me12b}&\\
  \hline
  $\Omega(n)$   & \cite{dj01} & нижняя, если есть изолированный ноль\\
  \hline
  $\Omega\left(\frac{nk}{\log n \cdot \log nk}\right)$ & \cite{kog87} & верхняя, почти все функции, $k\le 2^{n/2}$\\
  \hline
  \textbf{$\Omega\left(\frac{nk}{\log n + \log k}\right)$} & эта работа & нижняя, почти все функции, $k\le 2^{n/2}$\\
  \hline
\end{tabular}
    \caption{Основные известные оценки сложности булевых функций.}
    \label{table:review}
\end{table}

В области анализа сложности булевых функций отметим также ряд известных верхних оценок полученных в работах \cite{dj02,mu06}.

Основные известные оценки сложности в подклассах приведены в таблице \ref{table:review}.

\section{Нижние оценки сложности}
Используем следующие понятия, определенные в работе \cite{me12b}:
\begin{ddef}
Ненулевой вектор $\alpha \in B_k$ назовем \textit{разложимым} по векторам $\alpha_1,\alpha_2, ...,\alpha_t$ если выполнены следующие равенства
\[
\alpha = \alpha_1 \vee \alpha_2 \vee ... \vee \alpha_t, \qquad <\!\bar{\alpha},\alpha_1\!> = <\!\bar{\alpha}, \alpha_2\!> = ... = <\!\bar{\alpha},\alpha_t\!> = 0,
\]
где под символом $<\!\alpha, \beta\!>$ понимается скалярное произведение векторов $\alpha$ и $\beta$.
\end{ddef}
\begin{ddef}
Ненулевой вектор $\alpha \in B_k$ назовем \textit{ортогонально разложимым} по векторам $\alpha_1, ...,\alpha_t$ если выполнены следующие равенства
\[\begin{cases}
\alpha = \alpha_1 \oplus \alpha_2 \oplus ... \oplus \alpha_t\\
\alpha = \alpha_1 \vee \alpha_2 \vee ... \vee \alpha_t
\end{cases}\]
\end{ddef}

Обозначим $\mathbb{I}_k$ вектор размерности $k$, состоящий только из единиц. В случае $\alpha = \mathbb{I}_k$ в условиях предыдущего определения назовем вектора $\{\alpha_i\}_{i=1}^t$ \textit{ортогональным разложением единицы}. Заметим, что если вектор $\alpha$ разложим по векторам $\{\alpha_i\}_{i\in I}$, то $\chi(\alpha) > \chi(\alpha_i), i\in I$.

Конъюнкция $K$ называется \textit{импликантой} функции $f(x_1,x_2,...,x_n)$, если $N_K \subseteq N_f$. Импликанта $K$ называется простой, если из $K$ не может быть вычеркнут ни один литерал так, чтобы полученная конъюнкция была импликантой $f(x_1,x_2,...,x_n)$. Термины \textit{простая импликанта} и \textit{несократимая конъюнкция} употребляются в работе как синонимы.

Обозначим $Z_f$ --- множество нулей функции $f$, а $N_f$ --- множество ее единиц. Определим вектора $\{e_i^k\}_{i=1}^{k}$ равенствами $e_{i+1}^k = \chi^{-1}(2^i)$, $0~\le~i~<~k$.

Пусть $E(t) = \{i:M_i^t = 1\}$, $Z(t) = [k]\setminus E(t)$. Отметим, что для приведенной функции $Z(t)\neq\varnothing$ и $E(t)\neq\varnothing$. Следуя работе \cite{dj01}, обозначим $\tilde{\theta}(i,j) = (\theta_1,\theta_2,...,\theta_n)$ точку булева куба такую, что $\theta_r = M_i^r$ при всех $r\in[n]\setminus \{j\}$ и $\theta_j = 1 - M_i^j$.

Далее будем полагать, что число единиц каждого столбца матрицы нулей всякой рассматриваемой далее функции не превосходит числа нулей. Отметим, что всякая булева функция может быть приведена к такому виду преобразованиями Шеннона-Поварова.

Для конъюнкции $K$ обозначим $rank^{+}\; K$ число положительных литералов конъюнкции; $rank^{-}\; K$ число отрицательных литералов; \[rank\;K = rank^{+}\;K + rank^{-}\;K.\]

В работе \cite{dj01} введено следующее определение
\begin{ddef}
Множество
\[\bigcup\limits_{i=1}^k\bigcup\limits_{j=1}^n \{\tilde\theta(i,j)\}\setminus \bar{N}_f\]
назовем \textit{множеством околонулевых точек}.
\end{ddef}

Обозначим указанное множество $\Theta$. Положим
\[\Theta^{1} = \bigcup\limits_{i=1}^k\bigcup_{\substack{M_i^j = 1\\ j\in \{1,2,...,n\}}} \{\tilde\theta(i,j)\}\setminus \bar{N}_f
\qquad \Theta^{0} = \bigcup\limits_{i=1}^k\bigcup_{\substack{M_i^j = 0\\ j\in \{1,2,...,n\}}} \{\tilde\theta(i,j)\}\setminus \bar{N}_f\]

А.~Г.~Дьяконовым в той же работе \cite{dj01} была доказана
\begin{lemma}
Для всех элементарных конъюнкций $K\in D^{\text{сокр}_f}$ ($D^{\text{сокр}_f}$ --- сокращенная ДНФ функции $f$) и всех $i \in [k]$ справедливо
\[\left|N_K \cap \bigcup_{j=1}^n \tilde\theta(i,j)\right| \le 1\]
\end{lemma}

\textit{Соседними} назовем булевы вектора одинаковой размерности, находящиеся на расстоянии 1 в метрике Хемминга. Рассмотрим произвольную приведенную булеву функцию $f(x_1,...,x_n)$, не имеющую соседних нулевых точек; обозначим $M_f$ ее матрицу нулей и зафиксируем ее произвольную ДНФ $D$.

Справедливы следующие утверждения, доказательства которых получены в работе \cite{me12a}

\begin{lemma}\label{literal_min_number_a}
Выделим в $D$ конъюнкции, содержащие литерал $x_i^{\sigma_i}$ в количестве $t$ штук. Литерал $x_i^{\sigma_i}$ входит еще как минимум в одну, отличную от выделенных, конъюнкцию из $D$, если при любом разрезании матрицы $M_f$ на $t$ частей $M_1,M_2,...,M_t$ существует такая матрица $M_j$, что для функции $\phi_j$ заданной матрицей $M_j$, ни одна из выделенных конъюнкций не определяет разбиение $x_i^{\sigma_i}$.
\end{lemma}
\begin{lemma}\label{literal_min_number_b}
Литерал $x_i^{\sigma_i}$ входит не менее чем в $t+1$ конъюнкцию каждой ДНФ функции $f$, если при любом разрезании матрицы $M_f$ на $t$ частей $M_1,M_2,...,M_t$, существует такая матрица $M_j$, что для функции $\phi_j$, заданной матрицей $M_j$, не существует набора литералов, образующего разбиение литерала $x_i^{\sigma_i}$.
\end{lemma}

Рассмотрим оценки сложности функций класс $P^n_k$. Далее через $\log$ обозначен натуральный логарифм.

\begin{theorema}
	Для почти всех функций класса $P_k^n$ ни одна из их дизъюнктивных нормальных форм состоящих только из простых импликант,
	не содержит импликант ранга $d$, если $d$ не лежит вне интервала $[\log k + O(\log\log k + \log \log n); \log nk - \Omega(\log\log nk)]$.
\end{theorema}
\begin{proof}
	Для доказательства теоремы достаточно заметить, что для вероятности $p$ существования простой импликанты ранга $d$ в ДНФ функции
	класса $P_k^n$ справедливо
	\[
	\begin{cases}
		p \le 2^d \dbinom nd \left(1 - 2^{-d}\right)^k\\
		p \le 2^d \dbinom nd \left(k/2^{d-1}\right)^d.
	\end{cases}
	\]
	Первое из неравенств системы следует из оценки вероятности для конъюнкции быть импликантой функции.
	Второе неравенство отвечает за простоту импликанты. В частности, ни при удалении любого литера из
	рассматриваемой конъюнкции она не должна более являться импликантой функции.
	
	Простые вычисления дают, что с вероятностью стремящейся к 1 выполнено
	\[
	p \ge \log k +O(\log \log n + \log\log k),
	\]
	и
	\[
	p \le \log nk - \Omega(\log\log nk),
	\]
	что завершает доказательство.
\end{proof}

\begin{theorema}
	Почти все функции класса $P_k^n$ не могут быть реализованы дизъюнктивынми нормальными формами состоящими меньше чем
	\[
		\frac{nk}{\log{nk}} (1+ O(\sqrt{\log nk}))
	\]
	конъюнкций
\end{theorema}
\begin{proof}
	Доказательство теоремы в целом основано на лемме ?? и предыдущей оценке. Согласно лемме, для заданой конъюнкции
	нам необходимо оценить число строк матрицы нулей функции, в которых подстроки соотвествующие конъюнкции, отличаются от ее
	сигнатуры не более чем на 1. В нашем случае это будет случайная величина удовлетворяющая условию ограниченных разностей и
	имеющее математическое ожидание не превосходящее ранга рассматриваемой конъюнкции. Согласно предыдущей теореме в сочетании
	с неравенством МакДиармида получаем оценку на число конъюнкций $N_k$ в виде
	\[
		N_k\ge\frac{nk}{\log{nk}} \left(1 + \Omega\left(\sqrt{\log nk}\right)\right),
	\]
	что завершает доказательство теоремы.
\end{proof}

\begin{theorema}
При $\log_n k = o(n)$ в классе функций $P_k^n$ существует функция, минимальная ДНФ реализация которой содержит не менее $N_k$ конъюнкций, где для $N_k$ выполнено
\[N_k \ge \Omega\left( \frac{nk \log n}{\log k}\right).\]
\end{theorema}
\begin{proof}
Указанную оценку можно получить оценив сложность покрытия околонулевых точек функции принимающей нулевые значения на одном из слоев булева куба.
\end{proof}
\section{Полу-эффект Шеннона и открытые вопросы}
Ключевым открытым вопросом в области сложности булевых функция является установление точных и асимптотически точных границ на поведение функций тех или иных классов.

Эффектом Шеннона состоит в том, что сложность реализации наиболее ``трудной'' булевой функции в заданном классе совпадает в асимптотике со сложностью реализации почти всех функций класса. В частности, указанная ситуация имеет место при исследовании сложности реализации булевых функций схемами без ограничений.

Полу-эффект Шеннона имеет место в случае, когда сложность реализации почти всех булевых функций имеет один и тот же порядок, но вообще говоря отличается от сложности реализации наиболее трудной функции. Такая ситуация имеет место при исследовании сложности реализации функций дизъюнктивными нормальными формами. Как несложно, заметить экстремальные оценки сложности достигаются на линейной функции, для реализации которой требуется $2^{n-1}$. В тоже время вариационный принцип устанавливает факт, состоящий в том, что почти все булевы функции имеют один порядок при реализации в ДНФ \cite{ng67}. Однако, вопрос какова эта сложность остается невыясненным до сих пор.

Отметим в этом направлении ряд работ, которые позволили достаточно точно его оценить \cite{gl64,ko83,pip03}.

В классе булевых функций с ограниченным числом нулей полу-эффект (и эффект) Шеннона реализуется при числе нулей не превосходящим $\log_2 n \phi(n)$, для некоторой положительной функции $\phi(n) \rightarrow +\infty$, при $n\rightarrow\infty$.

Из работы \cite{me12b} следует отсутствие эффекта Шеннона в окрестности числа нулей $k = \log n$, если в качестве меры сложности рассматривать число литералов входящих в ДНФ. При этом в оригинальной работе рассматривался, вообще говоря более узкий подкласс, чем $P_k^n$. Однако результаты без труда переносяться на более общий случай.

В рассматриваемых классах $P_k^n$ при достаточно большом $k$ эффекта Шеннона нет. Однако результаты этой работы фактически устанавливают, что порядок роста почти всех функций один, с точностью до константы. Тем не менее вопрос о полу-эффекте остается открытым.

Таким образом, авторам видится перенесение идей вариационного принципа Нигматуллина \cite{ng67} на анализа концентрации в классах $P^n_k$ одной из основных открытых задач в этой области.

\section{Заключение}
В работе получены новые нижние оценки на сложность дизъюнктивных нормальных форм булевых функций, обращающихся в ноль на относительно небольшом множестве из $k$ точек.

Для случая числа точек $k$ ограниченного полиномом от размерности пространства входных переменных, получен точный порядок роста сложности (числа конъюнкций, числа литералов) ДНФ булевых функций.

В заключительной части работы поставлен ряд открытых вопросов, касающихся реализации булевых функций в классе ДНФ, в частности вопрос о существовании полу-эффекта Шеннона в подклассах булевых функций с малым числом нулей.

\bibliographystyle{chicago}
\bibliography{shannon}

\end{document}